\documentclass{amsart}

\usepackage{xy,amssymb,diagrams}
\usepackage[hyperindex,pdfmark]{hyperref}

\xyoption{poly}
\xyoption{arc}
\xyoption{knot}
\xyoption{curve}
\xyoption{arrow}
\xyoption{all}

\makeindex

\newcommand{\dummy[1]}{{#1}}
\newcommand{\C}{\mathbb{C}}

\newcommand{\N}{\mathbb{N}}
\newcommand{\Z}{\mathbb{Z}}

\newcommand{\wis}[1]{{\text{\em \usefont{OT1}{cmtt}{m}{n} #1}}}

\newtheorem{theorem}{Theorem}[section]

\newtheorem{proposition}[theorem]{Proposition}

\theoremstyle{definition}

\theoremstyle{remark}

\begin{document}

%%%%%%%%%%%%%%%%%%%%%
%%%some commands
%%%%%%%%%%%%%%%%%%%%%
\sloppy
%%%%%%%%%%%%%%%%%%%%%
%%%end of commands
%%%%%%%%%%%%%%%%%%%%%

\title{Simple roots of deformed preprojective algebras}

% Remove or comment out any unused author tags.
\author{Lieven Le Bruyn}% first author
\address{Universiteit Antwerpen (UIA) \\ B-2610 Antwerp (Belgium)}
%\curraddr{}
\email{lebruyn@uia.ua.ac.be}
\urladdr{http://win-www.uia.ac.be/u/lebruyn/}

\maketitle

\begin{center} {\it For Idun Reiten on her 60th birthday.}
\end{center}

\begin{abstract}
In \cite{Crawley:1999} W. Crawley-Boevey gave a description of the set $\Sigma_{\lambda}$ consisting
of the dimension vectors of simple representations of the deformed preprojective algebra
$\Pi_{\lambda}$. In this note we present alternative descriptions of $\Sigma_{\lambda}$.
\end{abstract}

\section{Reduction to $\Pi_0$}

Recall that a quiver $\vec{Q}$ is a finite directed graph on a set of vertices $Q_v = \{ v_1,\hdots,
v_k \}$, having a finite set of arrows $Q_a = \{ a_1,\hdots,a_l \}$ where we allow both multiple
arrows between vertices and loops in vertices. The Euler form of $\vec{Q}$ is the bilinear form on
$\Z^k$ determined by the integral $k \times k$ matrix having as its $(i,j)$-entry
$\chi_{ij} = \delta_{ij} - \# \{ \text{arrows from $v_i$ to $v_j$} \}$. The double quiver $\bar{Q}$ of the
quiver $\vec{Q}$ is the quiver obtained by adjoining to every arrow $a \in Q_a$ an arrow $a^*$ in the
opposite direction. The path algebra $\C \bar{Q}$ has as $\C$-basis the set of all oriented paths
$p = a_{i_u} \hdots a_{i_1}$ of length $u \geq 1$ together with the vertex-idempotents
$e_i$ considered as paths of length zero. Multiplication in $\C \bar{Q}$ is induced
by concatenation (on the left) of paths. For rational numbers $\lambda_i$, the deformed
preprojective algebra is the quotient algebra
\[
\Pi_{\lambda} = \Pi_{\lambda}(\bar{Q}) = \dfrac{\C \bar{Q}}{(\sum_{a \in Q_a} [a,a^*] - \sum_{v_i \in Q_v} \lambda_i e_i)}
\]
The (difficult) problem of describing the set $\Sigma_{\lambda}$ of all dimension vectors of simple
representations of $\Pi_{\lambda}$ was solved by W. Crawley-Boevey in \cite{Crawley:1999}.
He proved that for $\alpha$ a positive root of $\bar{Q}$, $\alpha \in \Sigma_{\lambda}$ if and only if
\[
p(\alpha) > p(\beta_1) + \hdots + p(\beta_r) \]
for every decomposition $\alpha = \beta_1 + \hdots + \beta_r$ with $r \geq 2$ all all $\beta_i$
positive roots of $\bar{Q}$ such that $\lambda.\beta_i = 0$ and where $p(\beta) = 1-\chi(\beta,\beta)$.

 For a given dimension vector $\alpha = (a_1,\hdots,a_k) \in \N^k$ one defines the affine scheme
 $\wis{rep}_{\alpha}~\Pi_{\lambda}$ of $\alpha$-dimensional representations of $\Pi_{\lambda}$. There
 is a natural action of the basechange group $GL(\alpha) = \prod_{i=1}^k GL_{a_i}$ on this scheme
 and the corresponding quotient morphism
 \[
 \wis{rep}_{\alpha}~\Pi_{\lambda} \rOnto^{\pi} \wis{iss}_{\alpha}~\Pi_{\lambda} \]
sends a representation $V$ to the isomorphism class of the direct sum of its Jordan-H\"older factors.
 Let $\xi$ be a geometric point of $\wis{iss}_{\alpha}~\Pi_{\lambda}$,
then $\xi$ determines the isomorphism class of a semisimple $\alpha$-dimensional representation say with
decomposition
\[
M_{\xi} = S_1^{\oplus e_1} \oplus \hdots \oplus S_l^{\oplus e_l} \]
with the $S_i$ distinct simple representations of $\Pi_{\lambda}$ with dimension vector $\beta_i$ which
occurs in $M_{\xi}$ with multiplicity $e_i$. We say that $\xi$ is of representation type
$\tau = (e_1,\beta_1;\hdots;e_l,\beta_l)$. Construct a graph $G_B$ depending on the set of
simple dimension vectors $B = \{ \beta_1,\hdots,\beta_l \}$ having $l$ vertices
$\{ w_1,\hdots,x_l \}$ having $2p(\beta_i) = 2(1 - \chi(\beta_i,\beta_i)$ loops in vertex $w_i$ and
$- \chi(\beta_i,\beta_j)-\chi(\beta_j,\beta_i)$ edges between $w_i$ and $w_j$.

Let $\bar{Q}_{B}$ be the (double) quiver obtained from $G_{B}$
by replacing each solid edge by a pair of directed arrows with opposite ordering. In
\cite[\S 4]{Crawley:2001} W. Crawley-Boevey proved that there is an \'etale isomorphism between
a neighborhood of $\xi$ in $\wis{iss}_{\alpha}~\Pi_{\lambda}$ and a neighborhood of the
trivial representation $\overline{0}$ in $\wis{iss}_{\alpha_{\tau}}~\Pi_0(\bar{Q}_{B})$ where
$\alpha_{\tau} = (e_1,\hdots,e_l)$ determined by the multiplicities of the simple factors of
$M_{\xi}$.

The arguments in \cite[\S 4]{Crawley:2001} actually prove that there is a $GL(\alpha)$-equivariant
\'etale isomorphism between a neighborhood of the orbit of $M_{\xi}$ in $\wis{rep}_{\alpha}~\Pi_{\lambda}(\bar{Q})$
and a neighborhood of the orbit of $\overline{(1,0)}$ in the principal fiber bundle
\[
GL(\alpha) \times^{GL(\alpha_{\tau})} \wis{rep}_{\alpha_{\tau}}~\Pi_0(\bar{Q}_{\tau}) \]
Using the description of $\Sigma_{\lambda}$ it was proved in \cite{Crawley:1999} that
$\wis{iss}_{\alpha}~\Pi_{\lambda}$ is irreducible whenever $\alpha \in \Sigma$.

In this note we will give two alternative descriptions of the set $\Sigma_{\lambda}$ stressing the
fundamental role of the extended Dynkin quivers in the study of deformed preprojective
algebras. Both descriptions rely on the above irreducibility result so they do {\em not} give a short
proof of Crawley-Boevey's result unless an independent proof of irreducibility of $\wis{iss}_{\alpha}~\Pi_{\lambda}$
for all $\alpha \in \Sigma_{\lambda}$ is found. In the statement of the results we have therefore
separated the parts that depend on the irreducibility statement.

\begin{proposition} \label{reduce} Let $\xi$ be a geometric point of $\wis{iss}_{\alpha}~\Pi_{\lambda}$ of 
representation type $\tau = (e_1,\beta_1;\hdots;e_l,\beta_l)$. The following are equivalent
\begin{enumerate}
\item{Any neighborhood of $\xi$ in $\wis{iss}_{\alpha}~\Pi_{\lambda}$ contains a point of representation type $(1,\alpha)$ (whence, in
particular, $\alpha \in \Sigma_{\lambda}$).}
\item{$\alpha_{\tau} = (e_1,\hdots,e_l)$ is the dimension vector of a simple representation of
$\Pi_0(\bar{Q}_B)$.}
\item{Any neighborhood of $\overline{0}$ in $\wis{iss}_{\alpha_{\tau}}~\Pi_0(\bar{Q}_B)$ contains a
point of representation type $(1,\alpha_{\tau})$ (whence, in particular, $\alpha_{\tau}$ is the
dimension vector of a simple representation of $\Pi_0(\bar{Q}_B)$.}
\end{enumerate}
If moreover $\wis{iss}_{\alpha}~\Pi_{\lambda}$ is irreducible these statements are equivalent to
\begin{itemize}
\item{$\alpha \in \Sigma_{\lambda}$.}
\end{itemize}
\end{proposition}

\begin{proof}
By comparing the stabilizer subgroups of the closed orbits determined by corresponding points under
the \'etale isomorphism it follows that $(1) \Leftrightarrow (3)$ and clearly $(3) \Rightarrow (2)$.
Because the equations of $\Pi_0(\bar{Q}_B)$ are homogeneous there is a $\C^*$-action on
$\wis{rep}_{\alpha_{\tau}}~\Pi_0(\bar{Q}_B)$ (multiplying all matrices by $t \in \C^*$). The limit
point $t \rightarrow 0$ of any representation is the trivial representation. Starting from a simple
representation $V$, any neighborhood of $\overline{0}$ contains a point determined by $t.V$ for
suitable $t$ proving $(2) \Rightarrow (3)$. To prove that $\bullet \Rightarrow (1)$ observe that
the set of all points of representation type $(1,\alpha)$ form an open subset of $\wis{iss}_{\alpha}~\Pi_{\lambda}$
(follows from the \'etale local description), whence if $\wis{iss}_{\alpha}~\Pi_{\lambda}$ is irreducible
this set is dense.
\end{proof}

This result allows us to describe $\Sigma_{\lambda}$ inductively if we can determine the sets of simple
dimension vectors for preprojective algebras. The induction starts off by taking the positive roots $\alpha$
for $\vec{Q}$ minimal w.r.t. $\lambda.\alpha = 0$. It follows from the easier part of \cite{Crawley:1999}
that these $\alpha \in \Sigma_{\lambda}$.

\section{Genetic description of $\Sigma_0$}

In this section we start with the quiver $\vec{Q}$ and will give an inductive procedure to determine
$\Sigma_0$, the set of simple dimension vectors of $\Pi_0 = \Pi(\bar{Q})$.

Assume we have constructed a set $B = \{ \beta_1,\hdots,\beta_l \}$ with $\beta_i \in \Sigma_0$ (we
can take $\beta = \beta_i = \beta_j$ for $i \not= j$ provided $p(\beta) > 0$). We want to determine
the {\em minimal} linear combinations
\[
\alpha = e_1 \beta_1 + \hdots + e_l \beta_l \]
such that $\alpha \in \Sigma_0$. We will do this in terms of the graph $G_B$ constructed in the
previous section and the dimension vector $\alpha_{\tau} = (e_1,\hdots,e_l)$.

\begin{figure}
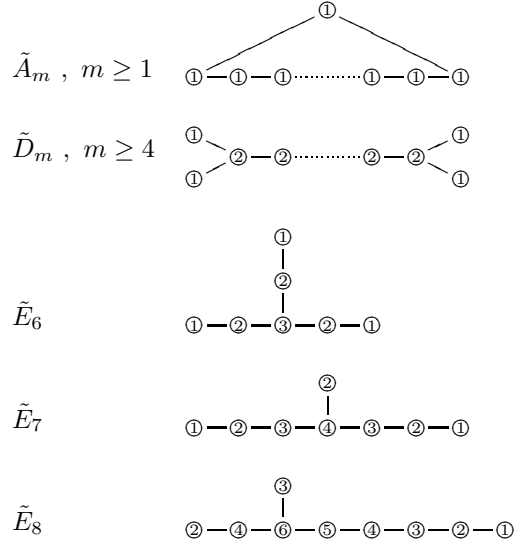

\[
\begin{array}{ll}
\tilde{A}_m~,~m \geq 1 &  \xy /r.07pc/:
\POS (0,0) *\cir<3pt>{}*+{\txt\tiny{1}} ="v1",
(60,30) *\cir<3pt>{}*+{\txt\tiny{1}} ="v0",
(20,0) *\cir<3pt>{}*+{\txt\tiny{1}} ="v2",
(40,0) *\cir<3pt>{}*+{\txt\tiny{1}} ="v3",
(80,0) *\cir<3pt>{}*+{\txt\tiny{1}} ="v4",
(100,0) *\cir<3pt>{}*+{\txt\tiny{1}} ="v5",
(120,0) *\cir<3pt>{}*+{\txt\tiny{1}} ="v6"
\POS"v1" \ar@{-} "v2"
\POS"v2" \ar@{-} "v3"
\POS"v3" \ar@{.} "v4"
\POS"v4" \ar@{-} "v5"
\POS"v5" \ar@{-} "v6"
\POS"v0" \ar@{-} "v1"
\POS"v0" \ar@{-} "v6"
\endxy
\\
& \\
\tilde{D}_m~,~m \geq 4 &  \xy /r.07pc/:
\POS (0,-10) *\cir<3pt>{}*+{\txt\tiny{1}} ="v1",
\POS (0,10) *\cir<3pt>{}*+{\txt\tiny{1}} ="v0",
(20,0) *\cir<3pt>{}*+{\txt\tiny{2}} ="v2",
(40,0) *\cir<3pt>{}*+{\txt\tiny{2}} ="v3",
(80,0) *\cir<3pt>{}*+{\txt\tiny{2}} ="v4",
(100,0) *\cir<3pt>{}*+{\txt\tiny{2}} ="v5",
(120,-10) *\cir<3pt>{}*+{\txt\tiny{1}} ="v6",
(120,10) *\cir<3pt>{}*+{\txt\tiny{1}} ="v7",
\POS"v1" \ar@{-} "v2"
\POS"v0" \ar@{-} "v2"
\POS"v2" \ar@{-} "v3"
\POS"v3" \ar@{.} "v4"
\POS"v4" \ar@{-} "v5"
\POS"v5" \ar@{-} "v6"
\POS"v5" \ar@{-} "v7"
\endxy \\
& \\
\tilde{E}_6 &  \xy /r.07pc/:
\POS (0,0) *\cir<3pt>{}*+{\txt\tiny{1}} ="v1",
(20,0) *\cir<3pt>{}*+{\txt\tiny{2}} ="v2",
(40,0) *\cir<3pt>{}*+{\txt\tiny{3}} ="v3",
(60,0) *\cir<3pt>{}*+{\txt\tiny{2}} ="v4",
(80,0) *\cir<3pt>{}*+{\txt\tiny{1}} ="v5",
(40,20) *\cir<3pt>{}*+{\txt\tiny{2}} ="v7",
(40,40) *\cir<3pt>{}*+{\txt\tiny{1}} ="v8"
\POS"v1" \ar@{-} "v2"
\POS"v2" \ar@{-} "v3"
\POS"v3" \ar@{-} "v4"
\POS"v4" \ar@{-} "v5"
\POS"v3" \ar@{-} "v7"
\POS"v7" \ar@{-} "v8"
\endxy \\
& \\
\tilde{E}_7 &  \xy /r.07pc/:
\POS (0,0) *\cir<3pt>{}*+{\txt\tiny{1}} ="v1",
(20,0) *\cir<3pt>{}*+{\txt\tiny{2}} ="v2",
(40,0) *\cir<3pt>{}*+{\txt\tiny{3}} ="v3",
(60,0) *\cir<3pt>{}*+{\txt\tiny{4}} ="v4",
(80,0) *\cir<3pt>{}*+{\txt\tiny{3}} ="v5",
(100,0) *\cir<3pt>{}*+{\txt\tiny{2}} ="v6",
(120,0) *\cir<3pt>{}*+{\txt\tiny{1}} ="v7",
(60,20) *\cir<3pt>{}*+{\txt\tiny{2}} ="v0"
\POS"v1" \ar@{-} "v2"
\POS"v2" \ar@{-} "v3"
\POS"v3" \ar@{-} "v4"
\POS"v4" \ar@{-} "v5"
\POS"v5" \ar@{-} "v6"
\POS"v6" \ar@{-} "v7"
\POS"v4" \ar@{-} "v0"
\endxy \\
& \\
\tilde{E}_8 &  \xy /r.07pc/:
\POS (0,0) *\cir<3pt>{}*+{\txt\tiny{2}} ="v1",
(20,0) *\cir<3pt>{}*+{\txt\tiny{4}} ="v2",
(40,0) *\cir<3pt>{}*+{\txt\tiny{6}} ="v3",
(60,0) *\cir<3pt>{}*+{\txt\tiny{5}} ="v4",
(80,0) *\cir<3pt>{}*+{\txt\tiny{4}} ="v5",
(100,0) *\cir<3pt>{}*+{\txt\tiny{3}} ="v6",
(120,0) *\cir<3pt>{}*+{\txt\tiny{2}} ="v8",
(140,0) *\cir<3pt>{}*+{\txt\tiny{1}} ="v9",
(40,20) *\cir<3pt>{}*+{\txt\tiny{3}} ="v0"
\POS"v1" \ar@{-} "v2"
\POS"v2" \ar@{-} "v3"
\POS"v3" \ar@{-} "v4"
\POS"v4" \ar@{-} "v5"
\POS"v5" \ar@{-} "v6"
\POS"v6" \ar@{-} "v8"
\POS"v3" \ar@{-} "v0"
\POS"v8" \ar@{-} "v9"
\endxy
\end{array}
\]
\caption{The tame settings.}
\label{tames}
\end{figure}

The tame settings are the couples $(D,\delta)$ where $D$ is an extended Dynkin diagram
and $\delta$ the corresponding imaginary root. The list of tame settings is given in 
figure~\ref{tames}.
We say that a tame setting $(D,\delta)$ is contained in $(G_B.\alpha_{\tau})$ if $D$ is a subgraph
of $G_B$ and if $\delta \leq \alpha_{\tau}$.

Recall from \cite{LeBruynProcesi:1990} that all polynomial invariants of quivers are generated by taking
traces along oriented cycles in the quiver. As a consequence, the coordinate algebra
$\C[\wis{iss}_{\alpha}~\Pi_0] = \C[\wis{rep}_{\alpha}~\Pi_0]^{GL(\alpha)}$ is
generated by traces in the quiver $\bar{Q}$. Note that non-trivial invariants exist whenever
$\alpha \in \Sigma_0$ and $\alpha$ is not a real root of $\vec{Q}$. The crucial ingredient in
our descriptions is the following technical result.

\begin{proposition} \label{tame} For $\alpha \in \Sigma_0$, if $\alpha$ is not a real root of $\vec{Q}$ and $\vec{Q}$ has only loops at vertices where
$\alpha$ is one,
then there is a non-loop tame setting $(\bar{D},\delta)$ contained in $(\bar{Q},\alpha)$.
\end{proposition}

\begin{proof}
Assume $\bar{Q}$ is a counterexample with a minimal number of vertices. There are 
at most two
directed arrows between two vertices ($(\tilde{A}_1,(1,1))$ is not contained) so we can define the
graph $G$ replacing a pair of directed arrows by a solid edge. Then, $G$ is a tree 
($(\tilde{A}_m,(1,\hdots,1))$ is not contained). 

We claim that the component of $\alpha$ for every internal (not a leaf) vertex is at least two. Assume
$v$ in internal and has dimension one, then any non-zero trace $tr(c)$ along a circuit in $\Gamma$
passing through $v$ (which must be the case by minimality of the counterexample) can be decomposed
as
\[
0 \not= tr(c) = tr(t_1)tr(t_2) \hdots tr(t_m) \]
where $t_i$ is part of the circuit along a subtree rooted at $v$. But then $tr(t_i) \not= 0$ when evaluated
at representations of the preprojective algebra of the corresponding subtree, contradicting
minimality of the counterexample.

Hence, $G$ is a binary tree ($(\tilde{D}_4,(2,1,1,1))$ is not contained) and even a star with at
most three arms ($(\tilde{D}_m,(2,\hdots,2,1,1,1,1))$ is not contained). If $G$ does not contain
$\tilde{E}_i$ for $6 \leq i \leq 8$ as subgraph, then $\bar{Q}$ is a Dynkin quiver and one knows
that in this case there are no nontrivial invariants, a contradiction.

If $\delta_v$ is the vertex-simple
concentrated in vertex $v$, we claim that
 \[
 \chi(\alpha,\delta_v)+\chi(\delta_v,\alpha) \leq 0 \]
 for every vertex $v$. Indeed, it follows from \cite{Crawley:1999b} that for any non-isomorphic
 simple $\Pi_0$-representations $V$ and $W$ of dimension vectors $\beta$ and $\gamma$ we have
 \[
 dim~Ext^1_{\Pi_0}(V,W) = - \chi(\beta,\gamma) - \chi(\gamma,\beta) \]
 Therefore,
twice the dimension of $\alpha$ at $v$ is smaller or equal to the sum of the dimensions of $\alpha$
in the two (maximum three) neighboring vertices.
Fill up the arm of $G$ corresponding to the longest arm
of $\tilde{E_i}$ with dimensions
starting with $1$ at the leaf  and proceeding by the rule that twice the dimension is equal to the sum
of the neighboring dimensions, then we obtain a dimension vector $\beta$ such that
\[
\delta_i \leq \beta \leq \alpha \]
where $\delta_i$ is the imaginary root of $\tilde{E}_i$, a contradiction.
\end{proof}

\begin{theorem} \label{genetic} With notations as above, we have
\begin{enumerate}
\item{$\alpha = e_1 \beta_1 + \hdots + e_l \beta_l \in \Sigma_0$ whenever $\delta=(e_1,\hdots,e_l)$
is the imaginary root of an extended Dynkin subgraph $D$ of $G_B$.}
\item{If moreover $\wis{iss}_{\alpha}~\Pi_0$ is irreducible for all
$\alpha \in \Sigma_0$, the
set $\Sigma_0$ is obtained by iterating the procedure in $(1)$ starting from the set of all real
roots of $\vec{Q}$.}
\end{enumerate}
\end{theorem}

\begin{proof}
(1) : There is a point $\xi \in \wis{iss}_{\alpha}~\Pi_0$ determined by a
semi-simple representation $M_{\xi}$ of representation type $\tau = (e_1,\beta_1;\hdots;e_l,\beta_l)$.
A neighborhood of $\xi$ is \'etale isomorphic to a neighborhood of $\overline{0}$ in
$\wis{iss}_{\delta}~\Pi_0(\bar{Q}_B)$. It is well known that $\wis{iss}_{\delta}~\Pi_0(\bar{D})$
contains points of representation type $(1,\delta)$ whence $\delta$ is a dimension vector of a
simple representation of $\Pi_0(\bar{Q}_B)$ (take a simple of $\Pi_0(\bar{D})$ and add zero matrices
for the remaining arrows). By proposition~\ref{reduce} it follows that $\alpha \in \Sigma_0$.

(2) : Let $\alpha \in \Sigma_0$ and take a decomposition (representation type)
\[
\alpha = d_1 \beta_1 + \hdots + d_l \beta_l \]
with all $\beta_i \in \Sigma_0$, $\beta_i < \alpha$ and $d = \sum_i d_i$ minimal.
Note that we can take all $d_i = 1$ whenever $p(\beta_i) > 0$ (as then there are
infinitely many non-isomorphic simples of dimension vector $\beta_i$). As a
consequence $G_B$ only has loops at vertices where $\alpha_{\tau}$ is equal to
one and $\alpha_{\tau}$ is a simple root for $\Pi_0(\bar{G_B})$ (here we
used irreducibility of $\wis{iss}_{\alpha}~\Pi_0$ in order to apply
proposition~\ref{reduce}. By proposition~\ref{tame} there is a non-loop tame subsetting
$(D,\delta)$ contained in $(G_B,\alpha_{\tau})$ and if $\delta = (e_1,\hdots,e_l)$
then we have a decomposition
\[
\alpha = (d_1-e_1) \beta_1 + \hdots + (d_l-e_l) \beta_l + 1.(\delta.\beta) \]
which has strictly smaller total number of multiplicities unless
$\alpha = \delta.\beta$. Induction on the total dimension finishes the proof.
\end{proof}

\section{Another description of $\Sigma_{\lambda}$}

In this section we reformulate the previous arguments in a more manageable
statement.

Take a non-trivial representation type $\tau = (d_1,\beta_1;\hdots;d_l,\beta_l)$ of $\alpha$ with all $\beta_i \in \Sigma_{\lambda}$. Let $\tau'$ be the 
representation type obtained from $\tau$ by replacing each $(d_i,\beta_i)$ by
$(1,\beta_i;\hdots;1,\beta_i)$ whenever $p(\beta_i) > 1$ (see the proof of
theorem~\ref{genetic}) and let $B'$ be the corresponding set os simple root
(some occurring more than once).

\begin{theorem} The following are equivalent
\begin{enumerate}
\item{$\alpha \in \Sigma_{\lambda}$ and $\wis{iss}_{\alpha}~\Pi_{\lambda}$ is
irreducible.}
\item{For all non-trivial representation types $\tau$ of $\alpha$ there is a
non-loop tame setting contained in $(G_{B'},\alpha_{\tau'})$.}
\end{enumerate}
\end{theorem}

\begin{proof}
$(2) \Rightarrow (1)$ : We claim that $(1,\alpha)$ is the unique maximal
representation type in the ordering of inclusion in Zariski-closures. Assume not and
let $\tau$ be another maximal type, then $\tau = \tau'$ and by proposition~\ref{tame}
there is a tame setting contained in $(G_B,\alpha_{\tau})$ but then there are
non-loop polynomial invariants, whence $\tau$ is not maximal.

$(1) \Rightarrow (2)$ : Follows from proposition~\ref{reduce} and proposition~\ref{tame}.
\end{proof}

Hence, the dimension vectors obtained from the genetic construction of theorem~\ref{genetic}
are exactly those $\alpha \in \Sigma_0$ such that $\wis{iss}_{\alpha}~\Pi_0$ is irreducible.

\par \vskip 5mm
\noindent
{\bf Acknowledgement : } I thank W. Crawley-Boevey for drawing my attention to
the circular argument used in the first version.

\end{document}